\documentclass{amsart}

\newcommand{\N}{\ensuremath{ \mathbf N }}

\newtheorem{theorem}{Theorem}
\newtheorem{lemma}{Lemma}
\newtheorem{corollary}{Corollary}
\newcommand{\bt}{\begin{theorem}}
\newcommand{\et}{\end{theorem}}
\newcommand{\bl}{\begin{lemma}}
\newcommand{\el}{\end{lemma}}
\newcommand{\bc}{\begin{corollary}}
\newcommand{\ec}{\end{corollary}}
\newcommand{\beq}{\begin{equation}}
\newcommand{\eeq}{\end{equation}}
\newcommand{\benum}{\begin{enumerate}}
\newcommand{\eenum}{\end{enumerate}}

\begin{document}

\title{Density of sets of natural numbers and the L{\' e}vy group}

\author{Melvyn B. Nathanson}
\address{Lehman College (CUNY),Bronx, New York 10468, and CUNY Graduate Center, New York, New York 10016}
\email{melvyn.nathanson@lehman.cuny.edu}

\author{Rohit Parikh}
\address{Brooklyn College (CUNY),
Brooklyn, New York 11210, and CUNY Graduate Center, New York, New York 10016}
\email{rparikh@gc.cuny.edu}

\thanks{The work of M.B.N. was supported in part by grants from the NSA Mathematical Sciences Program and the PSC-CUNY Research Award Program.  The work of R. P. was supported in part by a grant from the PSC-CUNY Research Award Program.}

\maketitle

\begin{abstract}
Let $\N$ denote the set of positive integers.  The asymptotic density of the set $A \subseteq \N$ is $d(A) = \lim_{n\rightarrow\infty} |A\cap [1,n]|/n$, if this limit exists.  Let $ \mathcal{AD}$ denote the set of all sets of positive integers that have asymptotic density, and let $S_{\N}$ denote the set of all permutations of the positive integers \N.  The group $\mathcal{L}^{\sharp}$ consists of all permutations $f \in S_{\N}$ such that $A \in \mathcal{AD}$ if and only if $f(A) \in \mathcal{AD}$, and  the group $\mathcal{L}^{\ast}$ consists of all permutations $f \in \mathcal{L}^{\sharp}$ such that $d(f(A)) = d(A)$ for all  $A \in \mathcal{AD}$.  Let $f:\N \rightarrow \N $ be a one-to-one function such that $d(f(\N))=1$ and, if $A \in \mathcal{AD}$, then $f(A) \in \mathcal{AD}$.  It is proved that $f$ must also preserve density, that is, $d(f(A)) = d(A)$ for all $A \in \mathcal{AD}$.  Thus, the groups $\mathcal{L}^{\sharp}$ and $\mathcal{L}^{\ast}$ coincide.
\end{abstract}

\section{Asymptotic density and permutations}
Let $A$ be a set of positive integers, and let 
\[
A(n) = \sum_{\stackrel{a\in A}{1 \leq a \leq n}} 1
\]
denote the counting function of the set $A$.
The {\em lower asymptotic density} of $A$ is
\[
d_L(A) = \liminf_{n \rightarrow \infty} \frac{A(n)}{n}.
\]
The {\em upper asymptotic density} of $A$ is
\[
d_U(A) = \limsup_{n \rightarrow \infty} \frac{A(n)}{n}.
\]
The set $A$ has {\em asymptotic density} $d(A)$
if the limit
\[
d(A) = \lim_{n \rightarrow \infty} \frac{A(n)}{n}
\]
exists.
The set $A$ has an asymptotic density if and only if $d_L(A) = d_U(A)$.
We denote by $\mathcal{AD}$  the set of all sets of positive integers that have asymptotic density, that is, 
\[
\mathcal{AD} = \{ A \subseteq \N : d_L(A) = d_U(A) \}.
\]

Let $S_{\N}$ denote the group of all permutations of the positive integers \N.
For any set $A \subseteq \N$ and permutation $g \in S_{\N}$, we let
\[
g(A) = \{g(a) : a \in A\}.
\]
Let $\mathcal{L}^{\ast}$ be the set of all permutations that preserve density, that is, $\mathcal{L}^{\ast}$ consists of all permutations $g\in S_{\N}$ such that
\benum
\item[(i)]
$A \in \mathcal{AD}$ if and only if $g(A) \in \mathcal{AD}$, and
\item[(ii)] 
$d(A) = d(g(A))$ for all $A \in \mathcal{AD}$.
\eenum
The set $\mathcal{L}^{\ast}$ is a subgroup of the infinite permutation group $S_{\N}$, and originated in work of Paul L{\' e}vy~\cite{levy51} in functional analysis.
This group and other related groups of permutations that preserve asymptotic density have been  investigated by Obata~\cite{obat88a,obat88b} and Bl{\"u}mlinger and Obata~\cite{blum-obat91}.

The L{\' e}vy group $\mathcal{L}^{\ast}$ is contained in the group $\mathcal{L}^{\sharp}$ that consists of all permutations $g\in S_{\N}$ such that
$A \in \mathcal{AD}$ if and only if  $g(A) \in \mathcal{AD}$, but that do not necessarily preserve the asymptotic density of every set $A \in \mathcal{AD}$.
The object of this note is to prove that $\mathcal{L}^{\ast} = \mathcal{L}^{\sharp}$.  Indeed, we prove the stronger result that if $f: \N\rightarrow \N$ is any one-to-one function, not necessarily a permutation, such that $A \in \mathcal{AD}$ implies that $f(A) \in \mathcal{AD}$, 
then also $d(f(A)) = \lambda d(A)$ for all $A \in \mathcal{AD}$, where $\lambda = d(f(\N))$.  In particular, if $f$ is a permutation, then $d(f(\N)) = d(\N) = 1$ and $d(f(A)) = d(A)$ for all $A \in \mathcal{AD}$.

\section{Permutations preserving density}

We beginning with the following ``intertwining lemma."

\bl     \label{levy:theorem:intertwining}
Let $A$ and $B$ be sets of integers such that $d(A) = d(B) = \gamma > 0.$  
Let $\{\varepsilon_k\}_{k=1}^{\infty}$ be a decreasing sequence of numbers such that $0 < \varepsilon_k < 1$ for all $k \geq 1$ and $\lim_{k\rightarrow \infty} \varepsilon_k = 0$.
Let  $\{M_k\}_{k=1}^{\infty}$ be a sequence of positive integers such that
\[
\left| \frac{A(n)}{n} - \gamma \right| < \varepsilon_k \text{ and } \left| \frac{B(n)}{n} - \gamma \right| < \varepsilon_k
\]
for all $n \geq M_k.$  
If $\{N_k\}_{k=1}^{\infty}$ is any sequence of integers satisfying
\[
M_{k-1} \leq N_{k-1} \leq \varepsilon_{k-1} N_k
\]
for all $k \geq 2$ and if
\[
C = \bigcup_{k=1}^{\infty}  (A \cap [N_{2k-1}+1, N_{2k}] ) \cup \bigcup_{k=1}^{\infty}  (B \cap [N_{2k}+1, N_{2k+1}] )
\] 
then 
\[
d(C)=\gamma.
\]
\el

\begin{proof}
If $N_k \leq m < n$, then
\[
(\gamma - \varepsilon_k)m < A(m) < (\gamma+\varepsilon_k)m
\]
\[
(\gamma - \varepsilon_k)n < A(n) < (\gamma+\varepsilon_k)n
\]
and so
\[
\gamma(n-m) - 2\varepsilon_k n < A(n)-A(m) < \gamma(n-m) + 2\varepsilon_k n.
\]
Similarly,
\[
\gamma(n-m) - 2\varepsilon_k n < B(n)-B(m) < \gamma(n-m) + 2\varepsilon_k n.
\]

Let $k \geq 2$ and $N_k < n \leq N_{k+1}$.  
If $k$ is odd, then 
\[
C \cap [1,n] = (A \cap [N_k+1, n]) \cup (B \cap [N_{k-1}+1, N_{k}]) \cup (C\cap [1,N_{k-1}])
\]
and so
\[
C(n) = A(n)-A(N_k) + B(N_k)-B(N_{k-1}) + C(N_{k-1}).
\]
If $k$ is even, then 
\[
C \cap [1,n] = (B \cap [N_k+1, n]) \cup (A \cap [N_{k-1}+1, N_{k}]) \cup (C\cap [1,N_{k-1}])
\]
and 
\[
C(n) = B(n)-B(N_k) + A(N_k)-A(N_{k-1}) + C(N_{k-1}).
\]
In both cases, since $N_{k-1} \leq \varepsilon_{k-1} N_k$, it follows that
\begin{align*}
C(n) & < \gamma(n-N_k) + 2\varepsilon_k n + \gamma(N_k - N_{k-1}) + 2\varepsilon_{k-1} N_k  + N_{k-1} \\
& < \gamma n + 5\varepsilon_{k-1} n
\end{align*}
and
\begin{align*}
C(n) & > \gamma(n-N_k) - 2\varepsilon_k n + \gamma(N_k - N_{k-1}) - 2\varepsilon_{k-1} N_k   \\
& > \gamma n - \gamma N_{k-1}  - 4\varepsilon_{k-1} n \\
& > \gamma n - 5\varepsilon_{k-1} n.
\end{align*}
Therefore,
\[
\left| \frac{C(n)}{n} - \gamma \right| < 5\varepsilon_{k-1}
\]
for all $n > N_k$, and so $d(C) = \gamma$.
\end{proof}

\bt        \label{levy:theorem:f-hat}
Let $f:\N\rightarrow \N$ be a one-to-one function such that if $A \in \mathcal{AD}$, then $f(A) \in \mathcal{AD}$, that is, if the set $A$ of positive integers has asymptotic density, then the set $f(A)$ also has asymptotic density.  Let $\lambda = d(f(\N))$.  If $\lambda = 0,$ then $d(f(A))=0$ for all $A \subseteq \N$.  If $\lambda > 0$, then there is a unique increasing function $\hat{f}:[0,1] \rightarrow [0,1]$ such that $\hat{f}(0) = 0$, $\hat{f}(1) = 1$, and 
\[
d(f(A)) = \lambda \hat{f}(d(A))
\]
for all $A \in \mathcal{AD}$.
\et

\begin{proof}
We shall prove that, for every set $A \in \mathcal{AD}$,  the asymptotic density of $f(A)$ depends only on the asymptotic density of $A$.  Equivalently, we shall prove that if $A,B \in \mathcal{AD}$ and $d(A) = d(B)$, then $d(f(A)) = d(f(B))$.

For $\gamma \in [0,1]$, let $A$ and $B$  be sets in $\mathcal{AD}$ such that $d(A) = d(B) = \gamma$.  Suppose that
\[
0 \leq d(f(A)) = \alpha <  \beta = d(f(B)) \leq 1.
\]
Let $\{\varepsilon_k\}_{k=1}^{\infty}$ be a decreasing sequence of numbers such that $0 < \varepsilon_k < 1$ for all $k \geq 1$ and $\lim_{k\rightarrow \infty} \varepsilon_k = 0$.  For every $k\geq 1$ there is a positive integer $M_k$ such that 
\begin{align*}
\left| \frac{A(n)}{n} - \gamma \right| & < \varepsilon_k \\
\left| \frac{B(n)}{n} - \gamma \right| & < \varepsilon_k \\
\left| \frac{f(A)(n)}{n} - \alpha \right| & < \varepsilon_k \\
\left| \frac{f(B)(n)}{n} - \beta \right| & < \varepsilon_k \\
\end{align*}
for all $n \geq M_k.$  
By Lemma~\ref{levy:theorem:intertwining},  if $\{N_k\}_{k=1}^{\infty}$ is any sequence of integers satisfying
\beq  \label{levy:MNineq}
M_{k-1} \leq N_{k-1} \leq \varepsilon_{k-1} N_k
\eeq
for all $k \geq 2$ and if
\beq \label{levy:setC}
C = \bigcup_{k=1}^{\infty}  (A \cap [N_{2k-1}+1, N_{2k}] ) \cup \bigcup_{k=1}^{\infty}  (B \cap [N_{2k}+1, N_{2k+1}] )
\eeq
then $d(C)=\gamma$. 

We shall construct a sequence $\{N_k\}_{k=1}^{\infty}$ satisfying~\eqref{levy:MNineq} such that the associated set $C$ satisfies $d(C) = \gamma$, but $d_L(f(C)) \leq \alpha$ and $d_U(f(C)) \geq \beta$.  This implies that the set $f(C)$ does not have asymptotic density, which is impossible since the function $f$ maps $\mathcal{AD}$ into $\mathcal{AD}$.  

The sequence $\{N_k\}_{k=1}^{\infty}$ and a related sequence $\{L_k\}_{k=1}^{\infty}$ will be constructed inductively.  We remark that since the function $f$ is one-to-one, it follows that for every positive integer $L$, there is an integer $N$ such that $f(n) \leq L$ only if $n \leq N$, and so
\[
f(C) \cap [1,L]  = f(C \cap [1,N]) \cap [1,L].
\]
Let $N_1 = L_1 = M_1$.
Let $k \geq 2$ and suppose that we have constructed sequences $N_1 < \cdots < N_{k-1}$ and $L_1 < \cdots < L_{k-1}$.  Choose an integer
\[
L_k > \max(L_{k-1},M_k)
\]
such that $\varepsilon_{k-1}  L_k >  N_{k-1}$.
By the remark, there exists an integer $N_k >  L_k$  such that $f(n) \leq L_k$ only if $n \leq N_k$.  Then 
\[
\varepsilon_{k-1} N_k  > \varepsilon_{k-1} L_k > N_{k-1}.
\] 
We use the sequence $\{N_k\}_{k=1}^{\infty}$ to construct the set $C$ according to formula~\eqref{levy:setC}.  

For $k \geq 1$ we have
\begin{align*}
f(C) \cap [1,L_{2k}] 
& = f(C \cap [1,N_{2k}]) \cap [1,L_{2k}] \\
& = \left(  \left( f(C \cap [1,N_{2k-1}]) \right) \cap [1,L_{2k}]\right)   \cup   \left( \left( f(A \cap [N_{2k-1}+1,N_{2k}]) \cap [1,L_{2k}] \right)\right) \\
& \subseteq f([1,N_{2k-1}]) \cup \left( f(A) \cap [1,L_{2k}] \right) 
\end{align*}
and so
\[
f(C)(L_{2k}) \leq  f(A)(L_{2k}) + N_{2k-1} .
\]
It follows that 
\[
\frac{f(C)(L_{2k})}{L_{2k}} \leq \frac{ f(A)(L_{2k})  + N_{2k-1}}{L_{2k}}  < \alpha + 2\varepsilon_{2k-1}.
\]
Therefore,
\[
d_L(f(C)) = \liminf_{n\rightarrow\infty} \frac{f(C)(n)}{n} \leq \liminf_{k\rightarrow\infty} \frac{f(C)(L_{2k})}{L_{2k}} \leq \alpha.
\]

Similarly,
\begin{align*}
f(C) \cap [1,L_{2k+1}] 
& \supseteq  f(B \cap [N_{2k}+1,N_{2k+1}]) \cap [1,L_{2k+1}]  \\
& = \left( f(B \cap [1,N_{2k+1}]) \cap [1,L_{2k+1}] \right) \setminus \left( f(B \cap [1,N_{2k}]) \cap [1,L_{2k+1}] \right) \\
& \supseteq \left( f(B) \cap [1,L_{2k+1}] \right) \setminus  f([1,N_{2k}]) 
\end{align*}
and so
\[
f(C)(L_{2k+1}) \geq  f(B)(L_{2k+1}) - N_{2k}.
\]
It follows that 
\[
\frac{f(C)(L_{2k+1})}{L_{2k+1}} \geq \frac{ f(B)(L_{2k+1}) - N_{2k}}{L_{2k+1}} > \beta - 2\varepsilon_{2k}
\]
and so
\[
d_U(f(C)) = \limsup_{n\rightarrow\infty} \frac{f(C)(n)}{n} \geq \limsup_{k\rightarrow\infty} \frac{f(C)(L_{2k+1})}{L_{2k+1}}\geq \beta.
\]
The inequality
\[
d_L(f(C)) \leq  \alpha < \beta  \leq d_U(f(C))
\]
contradicts the fact that $f(C)$ has  asymptotic density, and so $d(f(A)) = d(f(B))$.

If $\lambda = d(f(\N)) = 0$, then $d(f(A))=0$ for every set $A \subseteq \N$.  Suppose that $\lambda > 0$.  Define the function $\hat{f}$ by
\[
\hat{f}(\alpha) = \frac{d(f(A))}{\lambda}
\]
where $A \subseteq \N$ and $d(A) = \alpha$.  This is well-defined, since $d(f(A)) = d(f(A'))$ if $d(A) = d(A')$.  
Let $0 \leq \alpha \leq \beta \leq 1$.  There exist sets $A \subseteq B \subseteq \N$ such that $d(A) = \alpha$ and $d(B) = \beta$.  Since $f(A) \subseteq f(B) \subseteq f(\N)$, it follows that
\[
0 \leq d(f(A)) \leq d(f(B)) \leq d(f(\N)) = \lambda
\]
and so
\[
0 \leq  \hat{f}(\alpha)\leq \hat{f}(\beta)  \leq 1.
\]
Thus, $\hat{f}: [0,1]\rightarrow [0,1]$ is an increasing function with $\hat{f}(0) = d(f(\emptyset))= 0$ and $\hat{f}(1) = d(f(\N))/\lambda = 1$.  This completes the proof.
\end{proof}

\bt
Let $f:\N\rightarrow \N$ be a one-to-one function such that if the set $A$ of positive integers has asymptotic density, then the set $f(A)$ also has asymptotic density.  Let $\lambda = d(f(\N))$.  Then
\[
d(f(A)) = \lambda d(A)
\]
for all $A \in \mathcal{AD}$.
\et

\begin{proof}
If $\lambda = 0$, then $d(f(A)) = 0$ for all $A \in \mathcal{AD}$ and the theorem is true.  

Suppose that $\lambda > 0$.  By Theorem~\ref{levy:theorem:f-hat}, there is an increasing function $\hat{f}: [0,1] \rightarrow [0,1]$ such that $d(f(A)) = \lambda \hat{f}(d(A))$ for all $A \in \mathcal{AD}$.  We shall prove that $\hat{f}(\alpha) = \alpha$ for all $\alpha \in [0,1]$.  Since $\hat{f}$ is increasing, it suffices to show that $\hat{f}(\alpha) = \alpha$ for all positive rational numbers $\alpha \in (0,1]$.  

Let $\alpha = r/s$, where $1 \leq r \leq s$.  For $i = 1,\ldots, s$, let $A_i = \{a\in \N : a \equiv i \pmod{s} \}$.  Let $A = \cup_{i=1}^{r} A_i$.  Then $d(A_i) = 1/s$ for $i=1,\ldots,s$ and $d(A) = r/s$.  Since the function $f$ is one-to-one, the set $f(A)$ is the disjoint union of the $r$ sets $f(A_1),\ldots, f(A_{r})$.  Similarly,  $f(\N)$ is the disjoint union of the $s$ sets $ f(A_1),\ldots, f(A_{s})$.   Since $A, A_1,\ldots, A_s \in \mathcal{AD}$, it follows that $f(A), f(A_1),\ldots, f(A_s) \in \mathcal{AD}$, and
\[
\lambda = d(f(\N)) = \sum_{i=1}^s d(f(A_i))  =  \lambda s\hat{f}(1/s).
\]
Then
\[
\hat{f}(1/s) = \frac{1 }{s}
\]
and
\[
\hat{f}(\alpha) = \frac{d(f(A))}{\lambda} = \frac{1}{\lambda}\sum_{i=1}^r d(f(A_i)) =   \sum_{i=1}^r \hat{f}(d(A_i)) =  r\hat{f}(1/s) = \frac{r}{s} = \alpha.
\]
This completes the proof.
\end{proof}

\emph{Remark.}  The Levy group $\mathcal{L}^{\sharp}$ consists of all permutations $f \in S_{\N}$ such that $A \in\mathcal{AD}$ if and only if $f(A) \in\mathcal{AD}$.   We can also consider the semigroup $\mathcal{S}^{\sharp}$ consisting of all permutations $f \in S_{\N}$ such that $A \in\mathcal{AD}$ implies $f(A) \in\mathcal{AD}$.  The group $\mathcal{L}^{\sharp}$ is a subsemigroup of $\mathcal{S}^{\sharp}$.  It is natural to ask if $\mathcal{L}^{\sharp} = \mathcal{S}^{\sharp}$.  Equivalently, if $A$ is a set of positive integers such that $f(A)$ has asymptotic density for some $f \in \mathcal{S}^{\sharp}$, then does $A$ have asymptotic density?

\emph{Acknowledgements.} The authors thank Roman Kuznets, Brooke Orosz, and Samer Salame for many useful discussions.

\providecommand{\bysame}{\leavevmode\hbox to3em{\hrulefill}\thinspace}
\providecommand{\MR}{\relax\ifhmode\unskip\space\fi MR }
\providecommand{\MRhref}[2]{%
  \href{http://www.ams.org/mathscinet-getitem?mr=#1}{#2}
}
\providecommand{\href}[2]{#2}


\begin{thebibliography}{1}

\bibitem{blum-obat91}
M.~Bl{\"u}mlinger and N.~Obata, \emph{Permutations preserving {C}es\`aro mean,
  densities of natural numbers and uniform distribution of sequences}, Ann.
  Inst. Fourier (Grenoble) \textbf{41} (1991), no.~3, 665--678.

\bibitem{levy51}
Paul L{\'e}vy, \emph{Probl\`emes concrets d'analyse fonctionnelle. {A}vec un
  compl\'ement sur les fonctionnelles analytiques par {F}. {P}ellegrino},
  Gauthier-Villars, Paris, 1951, 2d ed.

\bibitem{obat88b}
Nobuaki Obata, \emph{Density of natural numbers and the {L}\'evy group}, J.
  Number Theory \textbf{30} (1988), no.~3, 288--297.

\bibitem{obat88a}
\bysame, \emph{A note on certain permutation groups in the infinite-dimensional
  rotation group}, Nagoya Math. J. \textbf{109} (1988), 91--107.

\end{thebibliography}
\end{document}